\documentclass[12pt,a4paper,twoside,reqno]{amsart}

\usepackage{amsfonts}
\usepackage{amssymb}
\usepackage{amsmath}    
\usepackage{amsthm}     
\usepackage{amscd}      
\usepackage{euscript}
\usepackage[matrix,arrow,curve]{xy}

{
       \newtheorem{theorem}{Theorem}[section]
       \newtheorem{proposition}[theorem]{Proposition}
       \newtheorem{lemma}[theorem]{Lemma}

      \newtheorem{corollary}[theorem]{Corollary}

{\theoremstyle{definition}

       \newtheorem{remark}[theorem]{Remark}
}

\newcommand{\RR}{{\mathbb{R}}}

\newcommand{\CC}{{\mathbb{C}}}

\newcommand{\PP}{{\mathbb{C} \mathbb{P}}}
\newcommand{\ZZ}{{\mathbb{Z}}}

\newcommand{\cX}{{\mathcal{X}}}

\newcommand{\p}{\partial}

\newcommand{\la}{\langle}
\newcommand{\ra}{\rangle}

\begin{document}

\title
[ Soliton-type metrics \& K\"ahler-Ricci flow
on sympl. quo.]
{Soliton-type metrics and K\"ahler-Ricci flow \\
on symplectic quotients}

\author{Gabriele La Nave and Gang Tian }
\maketitle
\tableofcontents

\section{Introduction}

In this paper, we first propose an interpretation of
the K\"ahler-Ricci flow on a manifold $X$ as an exact elliptic
equation of Einstein type on a manifold $M$ of which $X$ is one of
the (K\"ahler) symplectic reductions via a (non-trivial) torus action. There are plenty of such manifolds
(e.g. any line bundle on $X$ will do).

More precisely, let $M$ be a compact K\"ahler manifold which admits a Hamiltonian
$S^1$-action by holomorphic automorphisms and let $V$ be the vector field generating
such an action. Then there is a moment map $\mu: M\mapsto \sqrt{-1} \RR$ for this action.
Assume that $[0,\bar \tau]\subset \RR$
consists of regular values of $-\sqrt{-1} \mu$ and for $\tau \in [0,\bar\tau]$,
$X_\tau=\mu^{-1}(\sqrt{-1} \,\tau)/ S^1$ be the symplectic quotient of $M$ by this action.
All these $X_\tau$ are biholomorphic to each other. We consider K\"ahler metrics which are invariant under the $S^1$-action.
As usual, given a K\"ahler metric $g$, we denote by $\omega_g$ its K\"ahler form and ${\rm Ric}(g)$ its Ricci curvature form.

Our first result (cf. Theorem \ref{theorem:twistedricciflat}, \ref{theorem:converse} and Lemma \ref{meancurvature}) states, loosely speaking,
that the normalized K\"ahler-Ricci flow $\p _t \omega_g =-{\rm Ric} +\lambda \omega_g$ on $X_\tau$,
where $\lambda $ is a constant and $\tau \in \RR$, is equivalent to the system of
equations on $M$:
\begin{eqnarray}
\label{eq:twist-einstein}
\left\{\begin{aligned}& {\rm Ric}(g) +\frac{\sqrt{-1}}{2} \partial\overline{\partial} \left (\log (|V|_g^2)+f\right) =\lambda \,\omega_g\\
&\frac{d\tau}{dt} = -\frac{H(\tau)}{4|V|_g}+\frac{|V|_g^2}{4} \frac{\p f}{\p \tau}\end{aligned}\right.
\end{eqnarray}
for some function $f$ such that $f=f\cdot(-\sqrt{-1} \mu)$, that is, $f$ depends only on $\tau$, satisfying
$$-R(h)+n\lambda -\frac {\p f}{\p \tau}<0.$$
Here $J$ denotes the complex structure on $M$ and $H(\tau)$ denotes the mean curvature of the hypersurface $Y_\tau :=
\mu ^{-1}(\sqrt{-1} \tau) \subset M$ with respect to the metric $g$, which we
require to be $S^1$-invariant. Also we note that $R(g)$ is the scalar curvature of $g$.
We will call $g$ a {\bf V-soliton metric} if it is K\"ahler and satisfies:
\begin{equation}
\label{eq:twist-einstein2}
{\rm Ric} (g) + \frac{\sqrt{-1}}{2} \partial\overline{\partial} \left( \log (|V|_g^2)+f\right) =\lambda \,\omega_g.
\end{equation}

Such a $V$-soliton metric can be regarded as a generalization of K\"ahler-Einstein metrics or K\"ahler-Ricci solitons. Similarly to the case of K\"ahler-Einstein metrics, we can reduce (\ref{eq:twist-einstein2}) to a scalar equation on K\"ahler potentials, which is of Monge-Ampere type. To be more explicit, we fix a K\"ahler metric $g_0$ with K\"ahler form $\omega_0$ and write
$$\omega_g = \omega_{0} + \frac{\sqrt{-1}}{2} \p\bar\p u.$$
We can prove that if $M$ is compact\footnote{Even if $M$ is non-compact, the same conclusion still holds for the solutions of \eqref{eq:twist-einstein2}
with appropriate asymptotic behaviors.}, then the above $V$-soliton equation is equivalent to
the following scalar equation on $u$:
\begin{equation}
\label{eq:v-twisted-MA}
(\omega_{0}+ \frac{\sqrt{-1}}{2} \p \bar \p u )^n= |V|_g ^2 e ^{F+f-\lambda u}\omega_0^n,
\end{equation}
where $F$ is determined by $Ric (g) - \lambda\, g=\frac{1}{2} d (J d F)$.  We call the equation (\ref{eq:v-twisted-MA})
{\bf scalar $V$-soliton equation}. This equation is of complex Monge-Ampere type.

In the second part of this paper, we prove some preliminary results towards establishing existence of solutions
for (\ref{eq:v-twisted-MA}) on a compact K\"ahler manifold $M$.
In a forthcoming paper \cite{gatian}, we will establish an existence theorem for (\ref{eq:v-twisted-MA}).

This interpretation can be also extended to any symplectic quotients by more general groups.
An holomorphic Hamiltonian action of a Lie group $G$ on a manifold $M$ comes with a moment
map $\mu : M \to \mathfrak g ^*$: For every coadjoint orbit in the dual
Lie algebra of $G$, $\tau \subset\mathfrak g ^*$, there is a K\"ahler quotient $X_\tau := \mu ^{-1}
(\tau)/G$, as above, we can have an elliptic equation on $M$ whose solutions can descend to
solutions of the K\"ahler-Ricci flow on $X_\tau$ (cf. Theorem \ref{theorem:twistedricciflat}).

Our work was inspired by Perelman's groundbreaking work on Ricci flow. He
gave a formal interpretation of the (backwards) Ricci
flow on a manifold $M$ in terms of an {\it asymptotically} Ricci flat metric on $M\times
S^N\times \RR _+$, namely, the warped product metric $\tilde g=
g(\tau) + \tau g_{S^N} + \left( \frac{N}{2\tau} +R\right) d\tau ^2$,
where $g(\tau)$ solves the backward Ricci flow and $g_{S^N}$ is the
metric on $S^N$ with constant curvature equal to $\frac{1}{2N}$.
This allows him to heuristically interpret the
monotonicity of the reduced volume purely in terms of an analogue of Bishop-Gromov's
volume comparison theorem for asymptotically Ricci flat metrics. One of the major hurdles for
turning his heuristic description into a powerful tool to study the Ricci flow
is that the metric on $M\times
S^N\times \RR _+$ is only asymptotically Ricci-flat in $N$ (the dimension of the sphere).
Here we give a precise interpretation in the case of the K\"aher-Ricci flow.

One of our main motivations for this interpretation is to study singularity formation of
the K\"ahler-Ricci flow on a manifold with indefinite $c_1(M)$. A singularity can occur when the manifold is forced by the
flow to undertake a birational transformation. A large class of birational transformations can be constructed through
symplectic quotients. Then our interpretation may reduce studying singularity of the K\"ahler-Ricci flow on quotients
to studying an elliptic problem on $M$ which should be easier. In a subsequent paper,
we will discuss how our method can be applied to the study singularity formation along the K\"ahler-Ricci flow and we will first
illustrate it in a concrete example (cf. section \ref{subsection:finitetime}). In fact, our method should be applicable to more general situations than it actually seems at
first sight, since there is an associated GIT quotient description for any given flip.

\section{GIT versus symplectic quotients}

\subsection{GIT quotients}\label{git}

>From GIT (Geometric Invariant Theory) (cf. \cite{mgit}), we know that
given an action of a \textit{reductive} group $G$ on a projective
manifold $M$ with polarization $L$, one can define the GIT quotient of
$M$ via $G$ by considering the Zariski open set $M^{ss}(L)\subset M$
consisting of semistable points of the action, on which $G$ still acts
and in fact one can take the quotient $M^{ss}/G$. This depends only on
the choice of the polarization $L$, and is denoted by $M//G$. There is a well-defined
holomorphic map $\pi: M^{ss}\mapsto M//G$.

Clearly changing $L$ changes the quotient just in a birational manner, and the
change in the GIT quotient as $L$ varies is very well understood
(cf. \cite{tha}). It was shown there that the way the varieties change is by means of birational
transformations called {\it flips}.

\subsection{Symplectic quotients}\label{symquotient}
Let $(M, \omega)$ be a symplectic manifold. Assume there is a group $G$
acting \textit{symplectically} on $M$. Then there exists a moment map:
$$\mu : M\to \mathfrak g * , $$
\noindent
where $\mathfrak g $ denotes the Lie algebra of $G$. This is
described as follows: if $W \in \mathfrak g$, then $\langle \mu(x), W \rangle$
is the Hamiltonian function which generates the
flow given by the action of $W$ on $M$.

In these circumstances one can perform the symplectic quotient,
defined as $X_\lambda:=\mu^{-1} (\lambda)/G$ for some $\lambda$ a $G$-orbit
in $\mathfrak g ^*$. This quotient is in fact a (smooth) symplectic
manifold if $\lambda$ is a regular value for $\mu$, by the
Marsden-Weinstein reduction theorem.

If $M$ is a K\"ahler manifold
associated to a (quasi-)projective variety, these quotients coincide
with the GIT quotients encountered earlier
(cf. \cite{kw}). Furthermore, by the convexity Theorem (cf. \cite{gs1}),
the image of the moment mapping is \textit{convex} and there is
therefore a chamber subdivision according to the critical points of
the moment map. When one passes the walls of this chamber subdivision,
the symplectic manifold undergoes a symplectic surgery akin to the
blowing-up in Algebraic Geometry, as proven by Guillemin and Sternberg
in \cite{gs3}.

\subsection{K\"ahler quotients and their variations}\label{kahlerquotients}

Let $G$ be a compact connected Lie group acting symplectically on $(M,g)$
via holomorphic isometries,
and let $\mu :  M \to {\mathfrak g}^*$ denote the moment map, where
$\mathfrak g = Lie (G)$. Denote
by $G_c$ the complexification of $G$.

We can think of $\mathfrak g$ as a sub-bundle of $TM$ (in fact of
$T\mu ^{-1}(\tau)$).
Let $Q_p(\tau) \subset T_p \mu ^{-1} (\tau)$ be the orthogonal complement (with
respect to the given K\"ahler metric $g$) of $\mathfrak g$. Hence
$T_pM= Q_p(\tau) \oplus {\mathfrak g}_p \oplus J{\mathfrak g}_p$ is an
orthogonal decomposition. One readily checks that $Q (\tau)$ is $J$-invariant
and that $\{ Q_p\}_{p\in \mu ^{-1}(\tau)}$ is a $G$-invariant
distribution.
Also observe that if $\pi _{\tau} : \mu^{-1}(\tau) \to X_{\tau}$ is the natural
projection, then $d\pi _{\tau} : Q (\tau) \to T X_{\tau}$ induces an isomorphism.

Recall that the complex structure $J_\tau$ on the K\"ahler reduction
$X_\tau$ is defined by the condition that $d\pi _\tau \circ J=
J_\tau \circ d\pi _\tau$ where $\pi _\tau: Y_\tau:=\mu ^{-1}(\tau) \to X_\tau$ is the natural projection.

The following lemma is well-known.

\begin{lemma}
Given any regular value $\tau$ there is a direct sum decomposition of $Q(\tau
 )\otimes \CC = Q(\tau )^{(1,0)}\oplus  Q(\tau )^{(0,1)}$ into the $+\sqrt{-1}$ and
$-\sqrt{-1}$-eigenspaces respectively. Then $d \pi _{\tau}$ induces an
isomorphism: $Q^{(1,0)}(\tau)\to T ^{(1,0)}X_{\tau }$. Moreover,
The induced complex structure $J_{\tau}$ on $X_{\tau}$ is integrable.
\end{lemma}

By imposing that $\pi
_{\tau} : Y_\tau \to X_{\tau}$ be a Riemannian submersion,
we can define a natural Riemannian metric $g_\tau$ on $X_{\tau}$. Note that
$$g_{{\tau}} (d\pi _{\tau} (W_1), d\pi _{\tau} (W_2)) = g (W_1,W_2),~~~\forall ~W_1, W_2 \in Q(\tau).$$
The metric $g_\tau$ is in fact Hermitian with respect to
$J_{\tau}$. If we denote by $ i _{\tau}:\mu ^{-1}(\tau )\to M$
the natural inclusion, we have:

\begin{lemma}
The metric $g_{{\tau}}$ on $X_{\tau}$ is K\"ahler and the
corresponding K\"ahler form $\omega _{\tau}$ satisfies:
$$\pi_\tau ^* \omega _{{\tau}} = i _{\tau} ^* \omega.$$
Moreover, if $G=S^1$, then for any interval $I\subset \sqrt{-1}\RR$ which contains $a$ and consists of regular values of $\mu$,
$\mu ^{-1} (I)$ is symplectically equivalent to $\mu ^{-1}(a) \times
I$ (at least in a neighborhood of $Y_a$) endowed with the symplectic structure $\pi_a ^*\omega _a +d((\tau-a)
\beta)$, where $\beta$ is a connection 1-form on the circle bundle $\pi_a :
\mu ^{-1}(a) \to X_a$. In particular, the reduced symplectic form on $X_{\tau}$ is equivalent to $\omega _a + (\tau-a) c_1$, where $c_1$ is the Chern class of the principal bundle $\mu^{-1}(a) \to X_a$.
\end{lemma}

\begin{proof}
By definition, $\omega _{{\tau}} (W,Z)= g_{{\tau}} (J_{{\tau}} W,Z)$. On the
other hand, if $\bar W$ and $\bar Z$ are the unique $G$-invariant
sections of $Q(\tau)$ such that $d\pi _{\tau } (\bar W)=W$ and $d\pi
_{\tau } (\bar Z)=Z$ respectively, then one has:
$$\begin{aligned}
\pi_\tau ^* \omega _{{\tau}} (\bar W, \bar Z) &=  g_{{\tau}}
(J_{{\tau}}d\pi _{\tau} (\bar  W) , d\pi _{\tau} (\bar Z))\circ \pi \\
& =g_{{\tau}}
(d\pi _{\tau}  (J \bar  W) , d\pi _{\tau} (\bar Z))\circ \pi_\tau \\& =g(J
\bar W, \bar Z) = i_{\tau } ^* \omega (\bar W, \bar Z )
\end{aligned}$$
it is now easy to see that if, say, $\bar W$ is not in $Q(\tau )$,
then both sides of the equation amount to zero.
This also shows closedness of $\omega _{{\tau}}$ since this identity
shows that $\pi_\tau ^* d \omega _{{\tau}} = i_{\tau} ^* d \omega =0$,
and the surjectivity of $\pi_\tau$ implies that $d\omega _{\tau } =0$.
The statement on the symplectic equivalence follows directly from the
uniqueness part of the coisotropic embedding theorem
(cf. \cite{w}), whereas the statement on the nature of $\omega _{X_{\tau}}$ is a mere consequence of the fact on $Y_\tau=\mu ^{-1}(\tau)$, the form $\pi_a ^*\omega _a +d((\tau-a) \beta)$ restricts to $\pi_a ^*\omega _a \mid _{Y_a} + (\tau-a) d\beta$, and clearly $d\beta=\pi _a^* c_1$ .
\end{proof}

This lemma is of course a special case of a theorem for
symplectic quotients (cf. \cite{gs3}). It is then natural (and essential for our
constructions to come) to ask oneself whether such a result carries
through to the complex structure of the K\"ahler quotients. This
turns out to be true (cf. \cite{kw}).

Specifically, one can prove that so long as the moment map does not
cross critical values, then the complex structure does not change. In
order to describe things a little more in depth, we need to introduce
some notation:

Let $\Phi _s :M\to M$ represent the gradient flow of the Morse
function $||\mu ||^2$ (where the norm is in the dual of the Lie
algebra $\mathfrak g ^*$), and set (cf \cite{kw}):
$$ M^{min} ( O) :=\left\{ x\in M : \; \; \lim _{s\to +\infty} \Phi_s (x)
  \cap \mu ^{-1} (O) \neq \emptyset   \right\}$$
and (cf.\cite{gs4})\begin{footnote} {In fact, it was proved by Kempf
  and Ness that $M^s$ is nothing other than the set of semistable
  points of the action of $G_c$ on $M$, thereby connecting the GIT
  quotient with the K\"ahler reduction.} \end{footnote}:
$$ M^{s} ( O) :=\left\{x\in M : \; \; G_c x
  \cap \mu ^{-1} (O) \neq \emptyset   \right\}$$
for any coadjoint orbit $O$. Then one can prove:
\begin{lemma}(cf. \cite{kw}, \cite{gs4})
$M^{min} ( O)$ and $M^s(0)$ are $G_c$-invariant complex submanifolds of $M$. Furthermore,
there are natural biholomorphisms between $M^{min} (O)/G_c$ and $\mu
^{-1} (O)/ G$ and between $M^s(O)/G_c$ and $\mu
^{-1} (O)/ G$.
\end{lemma}

For simplicity, we assume that $G=S^1$ and its Lie algebra is identified with $\RR$.\footnote{All the subsequent discussions go through for a general
$G$ which is a maximal compact subgroup of a complex linear group, such as $SL(N,\CC)$.} Then the moment map $\mu $ takes values in $\RR$ and
$$ M^{min} (t) :=\left\{ x\in M : \; \; \lim _{s\to +\infty} \Phi_s (x)
  \cap \mu ^{-1} (t) \neq \emptyset   \right\}.$$
Then $M^{min}(t)$ is acted upon by $\CC ^*$ and if $t$ is a regular value, the natural holomorphic projection: $M^{min} (t)\mapsto M^{min} (t)/\CC
^*$ descends to a biholomorphism between $M^{min} (t)/ \CC^*\simeq \mu ^{-1} (t)/S^1$. It follows
that the complex manifolds $\mu^{-1}(t_1)/ S^1$ and $\mu^{-1}(t_2)/ S^1$ are
biholomorphic to each other whenever $t_1$ and $t_2$ are in an interval which does not contain any critical values of
$\mu$. For the readers' convenience, we will give a direct proof of this fact.

\begin{proposition}\label{variationofcomplex}
If $V$ has no zeros in a neighborhood of $\mu^{-1}([a,a+t_0])$, then the 1-parameter group of diffeomorphisms $\phi _t :M
\to M$ generated by the vector field $U=\frac{JV}{|V|^2 _g}$ induces biholomorphisms
$\tilde \phi _t : X_a \to X_{t+a}$ for $t\in [0,t_0]$.
\end{proposition}

\begin{proof}
If we write $\phi(x,t)=\phi_t(x)$, then
\begin{equation}\label{onepar}
\left\{ \begin{aligned}& \frac{d\phi }{dt} = U ({\phi})\\
&\phi (x,0)=x
\end{aligned}\right.
\end{equation}
\noindent
Since $\nabla \mu = JV$, $U=\frac{\nabla \mu }{|\nabla \mu|^2}$ and consequently,
$\mu(\phi_t(x))=a+t$ whenever $\mu(x)=a$.

Clearly, through the natural projections $\pi _{t'}:
Y_{t'} \to X_{t'}$, where $Y_{t'}=\mu^{-1}(t')$, $\phi_t$ induce diffeomorphisms $\tilde \phi _t
:X_a\to X_{a+t}$.

We want to show that these diffeomorphisms are actually biholomorphic
maps. For this purpose, we need to show
$$d {\tilde \phi} _t (J _a \tilde Z) = J_{a+t} d {\tilde \phi} _t (\tilde Z)$$
for any vector field $\tilde Z$ of $X_a$.
Let $\psi _t: M\to M$ be an integral curve of the vector field $JV$. They are biholomorphic maps
and there is $\lambda : M\times \RR \to \RR$ such that $\phi (x, t) = \psi
 (x, \lambda (x,t))$, since the vector fields $U$ and $JT$ are parallel. In fact, $\lambda(x,t)$ satisfies
\begin{equation}\label{onepar}
 \left\{ \begin{aligned}& \frac{d\lambda }{dt} = \frac{1}{|V|^2}\\
&\lambda (x,0)=0.
\end{aligned}\right.
\end{equation}
It follows
$$d\phi _t (W )(p) =d \psi _{\lambda (t,p)} (W) + \Lambda (W) JV.$$
Define
$$r_t(x) = \phi(x,t-\mu(x)).$$
Clearly, this defines a retraction $r_{t} : \mu ^{-1} ([a,a+t_0]) \to \mu^{-1} (t)$.
In particular, $dr_t(W)=W$ for any $W\in T\mu ^{-1}(t)$. Here
we think of $T\mu ^{-1}(t)$ as a subspace of $TM$).
Also,  $JV$ lies in the kernel of $dr_t$.

Let $Z=\tilde Z^h\in T_pM$ be the horizontal lifting of $\tilde Z$, i.e.,
the unique vector perpendicular to $V$ and $JV$ such that $d\pi _a (Z)=\tilde Z$. Then we have

\begin{equation}
\begin{aligned}
d {\tilde \phi} _t (J _a \tilde Z)&=
d {\tilde \phi} _t
(d \pi _a (J
Z))\\
&= d \pi _{a+t} (d \phi _t (JZ)) \\
&=   d \pi _{a+t} (d \psi _{\lambda (p,t)}
(JZ) + \Lambda (JZ) JV)\\&=  d \pi _{a+t} \left(dr_{a+t} (d \psi _{\lambda (p,t)}
(JZ) + \Lambda (JZ) JV)\right)\\
&=  d \pi _{a+t}\left(dr_{a+t} \left(d \psi _{\lambda (p,t)}
(JZ)\right)\right)\\&= d \pi _{a+t} \left(dr_{a+t} (J d \psi _{\lambda (p,t)}(Z))\right).
\end{aligned}
\end{equation}
Here we used the facts that $\psi$ is holomorphic and $ker (d\pi _{a+t} )= \la V, JV\ra$.
On the other hand, we have

\begin{equation}
\begin{aligned}
J_{a+t} d\pi _{a+t} (d { \phi} _t (Z)) &=  d\pi _{a+t} ([J \;  d \phi _t
( Z)]^h)\\
&= d\pi _t \left([J \;\left(d \psi _{\lambda (t,p)} (Z) + \Lambda
  (Z) JV \right)]^h\right)\\&=
 d\pi _{a+t} d r_{a+t} \left([\left(d \psi _{\lambda (p,t)} (JZ) - \Lambda
  (Z) V \right)]^h\right)
  \\&=
 d\pi _{a+t} \left([\left(d \psi _{\lambda (p,t)} (JZ) - \Lambda
  (Z) V \right)]^h\right)\\&=
 d\pi _{a+t} \left(dr_{a+t}(d \psi _{\lambda (p,t)} (JZ))\right)
\end{aligned}
\end{equation}
since $dr_{a+t} (W)=W$ for any $W\in T\mu^{-1}(a+t)$. This completes the proof.
\end{proof}

As a consequence, we conclude that if all $a+t$, where $t\in [0,t_0]$, are regular values of $\mu$ and we denote by $F_{a+t}: M^{min}(a)\mapsto X_t$ the above induced
holomorphic map, then $F_{a+t}= \tilde \phi_{t}\cdot F_a$.

\section{$V$-soliton metrics and K\"ahler-Ricci flow on
  symplectic quotients}

\subsection{The Ricci flow on symplectic quotients}\label{ricciflowonsymplquotients}

We will start with the case of a complex torus $T_{\CC}=(\CC ^*)^N$
acting holomorphically on a
smooth K\"ahler manifold $M$ (in fact $M$ does not need to be
smooth, it could for instance be a K\"ahler space with {\it canonical singularities})

Note that $T_{\CC}$ is the complexification of the real torus $T=(S^1)^N$,
which acts on $M$ by Hamiltonian diffeomorphisms.
Let $\mu: M\to Lie(T)^* =\RR ^N$ be the associated moment map. We denote by $g$ an invariant K\"ahler metric on $M$.

Let $z_1,\cdots, z_n$ be holomorphic coordinates on the quotient
manifold $X_{a}$, and let $\tau _i$ be ``moment map coordinates'', i.e.,
$\mu=(\tau_1,\cdots, \tau_N)$, sometime, we simply identify $\mu$ with its value $\tau=(\tau_1,\cdots,\tau_N)$.
If $N=1$, write $\tau=\tau_1$. Clearly, we have $d\tau _k =i_{V_k} \omega$,
where $\{V_k\}$ is a basis of vector fields which generate the Hamiltonian action
of $T$ and correspond to an orthonormal basis of the Lie algebra of $T$. We can define
$1$-forms $\theta _1, \cdots, \theta _N$ by
$$\theta_k(V_l)=\delta_{kl}, ~\theta_k(J V_l)=0, ~\theta_k\big |_{Q}=0,$$
where $\nabla\tau_l$ denotes the gradient of $\tau_l$ with respect to $g$. By the definition of
the moment map, we have $\nabla\tau_l= J V_l$. In particular, $\nabla \tau_l$ is tangent to
orbits of the action by $T_{\CC}$.
\begin{lemma}
For the above local coordinates, we have $g(dz_i, d\tau_k)=0$, $g(dz_i,\theta_k)=0$ and $g(\theta_k,d\tau_l)=0$,
where $g$ also denotes the induced metric on the cotangent bundle of $M$.
\end{lemma}
\begin{proof}
We have
$$g(dz_i,d\tau_k)= dz_i(\nabla \tau_k) = (\nabla \tau_k)(z_i) =0.$$
Clearly, the second follows from the first since $J(dz_i) = \sqrt{-1} dz_i$.
For the third, we have
$$g(\theta_k, d\tau_l)= \theta_k(\nabla \tau_l) = \theta_k(JV_l)=0.$$
\end{proof}
It follows from this lemma and a direct computation that in the above local coordinates, we can write the K\"ahler metric $g$ on $M$ as:
\begin{equation}\label{invmetric}
g = h_{i\bar j} dz_i d\bar z _j+ w_{kl} d \tau _k d \tau
_l +
w^{kl} \theta _k \theta _l\end{equation}
where $w^{ij}= g(V_i, V_j)$ (this also shows that
the $w^{ij}$'s are globally defined) and $\{w^{ij}\}$ is a positive definite matrix
and $\{w_{ij}\}$ is its inverse. Also, in the above proof, we have used the fact that
$d\tau_k(\nabla \tau_i)= \omega(V_k,JV_i) = w^{ki}$.

Using $g(JV_i, W)=\omega _g(V_i, W)= d\tau _i(W)$, where $J$ denotes the complex structure of $M$,
we can deduce $-J\theta_i = w_{ij} d\tau _j $.  We can thus infer that $w_{ij}d\tau _j -
\sqrt{-1} \theta _i$ is of type $(1,0),$ and rewrite $g$ as
$$g= h_{i\bar j} dz_i d\bar z _j+w^{kl} (w_{ki}d\tau _i -
\sqrt{-1} \theta _k)( w_{lj}d\tau _j +\sqrt{-1} \theta _l).$$
Also, we have the decomposition: $T^{(1,0)}M = Q^{(1,0)} \oplus \langle w_{ij}d\tau _j -
\sqrt{-1} \theta _i \rangle$.

In the sequel we will need the following:

\begin{lemma}\label{connection}
One has:
\begin{equation}\label{closed}
d \theta _k = \sqrt{-1}\left\{ -\frac{1}{2} \frac{\partial {h_{i\bar
    j}}}{\partial \tau _k} dz_i\wedge d \bar z _j - \frac {\partial
  w_{ki}}{\partial z_j} d\tau _i\wedge dz_j + \frac {\partial
  w_{ki}}{\partial \bar z_j}d\tau _i\wedge d\bar z_j\right \}
\end{equation}
\end{lemma}
\begin{proof}
For simplicity, we will assume $N=1$ and write $\tau=\tau_1$, the proof for $N > 1$ is identical.
Observe that the K\"ahler form of $g$ is given by
$$\omega_g = \frac{\sqrt{-1}}{2} h_{i\bar j} dz_i \wedge d\bar z _j- d \tau \wedge \theta .$$
Since this is closed, we get
$$d\theta = -\frac{\sqrt{-1}}{2} \frac{\p h_{i\bar j}}{\p \tau} dz_i\wedge  d\bar z_j
+\beta,$$
where $\beta$ is a real $2$-form of the form
$$\beta=\sum _i\left( q_i dz_i\wedge d \tau + \bar q_i d\bar z_i\wedge d \tau\right) + r d \tau \wedge \theta .$$
Note that $r$ is a real function. Since the associated complex structure $J$ is integrable, the $(0,2)$-part of $d(w d\tau - \sqrt{-1} \theta)$
vanishes. This implies
$$0=[d(w d\tau - \sqrt{-1} \theta)]^{0,2} = [d w\wedge d\tau - \sqrt{-1}(q_i dz_i\wedge d \tau + \bar q_i d\bar z_i\wedge d \tau)]^{0,2},$$
consequently,
$$0=[(\frac{\partial w}{\partial \bar z_i } - \sqrt{-1} \bar q_i) d\bar z_i \wedge d\tau   ]^{0,2}.$$
Hence,
$$q_i = \sqrt{-1} \frac{\partial w}{\partial z_i }.$$
On the other hand, since $\theta$ is invariant under the action, $L_V\theta =0$, that is,
$d i_V\theta + i_V d \theta=0$. But $i_V \theta =1$, so $d\theta (V, JV)=0$, which implies that $r=0$
and consequently, the lemma is proved.
\end{proof}

We can now calculate the volume form in a {\it
holomorphic frame}, namely:
\begin{lemma}
There is a holomorphic frame for which the volume form of $\omega
_g$ equals to $\det(w)^{-1}\det(h)$.
\end{lemma}
\begin{proof}  We first show the following claim: There
exists (local) functions $f_{ik}$ such that $\gamma _k = f_{kl} dz_l +
w_{kl} d\tau _l - \sqrt{-1} \theta _k$ are holomorphic. This is clearly
equivalent to showing that there exist smooth functions $f_{ik}$ such
that $[d \gamma _i] ^{1,1}= 0$, i.e., $d \gamma _i$ is of type $(2,0)$.

Using the formula for $d\theta_i$ in the above lemma, we have
$$d(w_{ij}d\tau _j -\sqrt{-1} \theta _i)= \frac{\partial w_{ij}}{\partial \tau_k} d\tau_k\wedge d\tau_j + \beta^{(2,0)}_i
+ \beta^{(1,1)}_i,$$
where $\beta ^{(2,0)}_i$ and $\beta ^{(1,1)}_i$ are of type $(2,0)$ and $(1,1)$, respectively. Then
$$\beta ^{(1,1)}_i = - \frac{1}{2} \frac{\p h_{k\bar j}}{\p \tau _i} dz_k
\wedge d\bar z_j + 2 \frac{\p w_{ij}}{\p z _k} dz_k \wedge (d\tau _j +
\sqrt{-1} w^{jl} \theta _l).$$
Since each $w_{ij}d\tau _j -
\sqrt{-1} \theta _i$ is a $(1,0)$ form, $d (w_{ij}d\tau _j -
\sqrt{-1} \theta _i)$ has vanishing $(0,2)$-part. It follows that
$$[\frac{\partial w_{ij}}{\partial \tau_k} d\tau_k\wedge d\tau_j]^{0,2}=0.$$
Hence,
$$\frac{\partial w_{ij}}{\partial \tau_k}=\frac{\partial w_{ik}}{\partial \tau_j},$$
consequently, $[d \gamma _i] ^{1,1}=0$ if and only if
\begin{equation}
\label{eq:f-ik}
\left\{ \begin{aligned}
& \frac{\p f_{ik}}{\p \tau _j}-2 \frac{\p w_{ij}}{\p z_k}=0\\
& \frac{\p h_{k\bar j}}{\p \tau _i} + 2 \frac{\p f_{ik}}{\p \bar z _j}
=0\end{aligned}\right.\end{equation}
The integrability conditions for this system are
\begin{equation}
\frac{\p^2 w_{ij}}{\p z_k \p \tau_l}= \frac{\p^2 w_{il}}{\p z_k \p \tau_j}, ~~~~
\frac{\p^2 h_{i\bar j}}{\p \tau _k\p \bar z_l}=\frac{\p^2 h_{i\bar l}}{\p \tau _k\p \bar z_j},~~~~
4\,\frac{\p^2 w_{ij}}{\p z_k \p \bar z_l} = - \frac{\p^2 h_{k\bar l}}{\p \tau _i \p \tau _j}.
\end{equation}
The first identity follows easily from the above symmetry on $\frac{\partial w_{ij}}{\partial \tau_k}$. The second and third follow
from $d(d\theta_i)=0$ and the formula for $d\theta_i$. Hence, we can solve the equations in (\ref{eq:f-ik}) for $f_{ik}$ and our claim is proved.
We can therefore infere the existence of a local holomorphic frame $d z_1,\cdots, dz_n, \gamma_1, \cdots, \gamma_N$. In this local frame,
$\omega_g$ can be written as
\begin{equation}
\label{eq:rep-omega-again}
\begin{aligned}
&\frac{\sqrt{-1}}{2} h_{i\bar j} dz_i \wedge d\bar z _j -
d \tau _k \wedge \theta _k\\
 = &\frac{\sqrt{-1}}{2}\left( h_{i\bar j} dz_i
\wedge d\bar z _j + w^{ij} (w_{ik} d\tau _k -
\sqrt{-1} \theta _i) \wedge (w_{jl} d\tau _l +
\sqrt{-1} \theta _j)  \right)\\
=&\frac{\sqrt{-1}}{2}\left( (h_{i\bar j} + w^{kl} f_{ki}\bar f_{lj}) dz_i
\wedge d\bar z _j - w^{ij} (f_{ik} d z_k \wedge \bar \gamma_j + \bar f_{jl} \gamma_i \wedge d\bar z_l)
+ w^{ij} \gamma_i \wedge \bar \gamma_j\right )\\
\end{aligned}\nonumber
\end{equation}
It follows
$$\omega _g ^{n+N} =(n+N)! \left (\frac{\sqrt{-1}}{2}\right )^{n+N} \det(h) \det(w _{ij})^{-1}
\, dz \wedge d \bar z \wedge \gamma \wedge \bar \gamma.$$
The lemma is proved.
\end{proof}

\noindent
Next we compute the complex Hessian of any $T$-invariant function.
\begin{lemma}\label{complexhessian}
For any $T$-invariant function $\phi \in C^2(M)$, we have
\begin{equation}
\label{hessian}
\begin{aligned}
\p \bar \p \phi =&\sum\frac{\p^2\phi}{\p z_i\p\bar z_j} dz_i\wedge d\bar z_j
 +\frac{1}{4}\sum \frac{\p\phi}{\p \tau _l}w^{lk} \left( \frac{\partial {h_{i\bar
j}}}{\partial \tau _k} dz_i\wedge d \bar z _j \right)\\
-& \frac{1}{2}\sum \frac{\p}{\p z_k}
\left(\frac{\p \phi}{\p \tau _i} w^{ij}\right) dz_k\wedge (w_{jl}d\tau_l+\sqrt{-1} \theta _j)\\
+&\frac{1}{2}\sum \frac{\p} {\p \bar z_k}\left( \frac{\p \phi}{\p \tau _i}
  w^{ij}\right) d\bar z_k \wedge (w_{jl}d\tau_l-\sqrt{-1} \theta _j)\\
  +&\frac{\sqrt{-1}}{2}\sum \frac{\p}{\p \tau _k} \left( \frac{\p \phi}{\p \tau _i}
w^{ij}\right)d \tau _k \wedge \theta _j.
\end{aligned}
\end{equation}
\end{lemma}
\begin{proof} First
$$
d\phi = \sum \frac{\p \phi}{\p z_i} dz_i + \sum \frac{\p \phi}{\p
  \bar z_i} d\bar z_i + \sum \frac{\p \phi }{\p \tau _i} d\tau _ i$$
Then, using the fact that $J d\tau _i = w ^{ij} \theta _j$, we get
$$\begin{aligned}
d(J d \phi))
= &\frac{2}{\sqrt{-1}} \sum\frac{\p^2\phi}{\p z_i\p\bar z_j} dz_i\wedge d\bar z_j + \sum \left(d\left(\sum \frac{\p \phi}{\p \tau _i}
w^{ij}\right)\wedge \theta _j + \left(\frac{\p \phi}{\p \tau _i}
w^{ij} \right) d\theta _j\right)\\
+ &{\sqrt{-1}}\sum \left (\frac{\p^2\phi}{\p \tau_j\p z_i } d\tau_j\wedge dz_i
-\frac{\p^2\phi}{\p \tau_j\p \bar z_i } d\tau_j\wedge d\bar z_i\right) \end{aligned}$$
Then the lemma follows from the fact that $d (J d \phi) = -2 \sqrt{-1} \p\bar \p \phi$ and a direct computation (aided by
formula (\ref{closed})).
\end{proof}

We can obtain the following fact about Hamiltonian functions from the Lemma above:
\begin{corollary}\label{hamiltonianclass}
If $T=S^1$ and $\omega _g= \omega _{g_0} -\frac{1}{4}d (J d u)$ for some $S^1$-invariant
function $u$, and if $\mu: M \to \RR$ is a Hamiltonian function of the
$S^1$-action with respect to $\omega _{g_0}$, then $\tilde \mu:= \mu
+\frac{1}{4} w_0^{-1} \frac{\p u}{\p\tau}$ is a Hamiltonian function with respect to $\omega _g$, where
$w_0=g_0(V,V)$ and $V$ is the associated vector field of the $S^1$-action.
\end{corollary}
\begin{proof}
Indeed, using Lemma \ref{complexhessian}, we have
$$\begin{aligned}\omega _g (V,U)&= \omega _{g_0}(V,U) -\frac{1}{4} d ( J d u )(V,U)\\
&=d\mu (U) -\frac{1}{4}  \left( d(\frac{\p u}{\p\tau} \frac{1}{w_0})\wedge \theta \right)(V,U)\\
&= U(\mu+\frac{1}{4} \frac{\p u}{\p\tau} \frac{1}{w_0})\end{aligned}$$
\end{proof}

\noindent
If $X_{a}$ is smooth, one can compute the curvature of quotient metric $g_{a}$
in terms of $g$ on $M$ via the Gauss-Codazzi equations and then
uses O'Neill's formula (cf. \cite{on}) for the Riemannian submersion: $\mu
^{-1}(a) \to X_{a}$.
However, we shall perform our computations by exploring the
K\"ahlerian structures.

In order to prove the next theorem, we need the following:
\begin{lemma}\label{lemma:hamiltonianlaplacian}
The Hamiltonian functions $\tau _k$ satisfy
\begin{equation}\label{hamiltonianlaplacian}
w^{kl} \frac {\partial \log \det (h)} {\partial
  \tau _l}=\Delta _g \tau _k - \frac{\partial w^{kl}}{\partial \tau _l}
\end{equation}
where $\Delta_g$ denotes the Laplacian of $g$ on $M$. Note that the right side of the above is
independent of choices of coordinates $z_1,\cdots,z_n.$
\end{lemma}
\begin{proof} Taking trace of \eqref {hessian}, one gets
\begin{equation} \label{reallaplacian}\Delta _g f =  h^{i \bar j} \left( 4 \frac{\p ^2 f}{\p z_i \p \bar z_j} +
w^{kl} \frac{\p f}{\p \tau _k} \frac {\p h_{i\bar j}}  {\p \tau_l}\right) +
\frac{\p} {\p \tau _k}\left( \frac{\p f} {\p \tau _l}
w^{kl}\right).\end{equation}
It follows
$$\Delta _g \tau _k =  h^{i\bar j} \left( w^{kl}  \frac {\p h_{i\bar j}}  {\p \tau
_l}\right) + \frac{\p w^{kl}} {\p \tau _l}.$$
Then we get the desired identity by noticing that
$$\frac {\partial \log \det (h)}{\partial
\tau _l}=  h^{i j} \frac {\p h_{ij}}{\p \tau  _l}.$$
\end{proof}

Let $(M,g)$ be a K\"ahler manifold with a torus $T$-action by holomorphic isometries. Let $V_1,\cdots, V_N$ be a basis
of the Killing vector fields generating this $T$-action. As before, we can write the moment map in the form
$\mu=(\tau_1,\cdots,\tau_N)$, where $d\tau_k = i_{V_k} \omega_g$, and $w^{ij}= g(V_i,V_j)$.
Let $\phi_{\tau}: X_a\mapsto X_{a+\tau}$ be the biholomorphism defined in Proposition \ref{variationofcomplex}. Fix a unit
vector field $V_\tau = \sum_i b_i V_i$, then we get an one-parameter family of metrics on $X_a$: $h(\tau)=\phi^*_{\tau} g_{a+\tau}$
so long as there are no critical points of $\mu$ in $\{ a + s b~|~ 0\le s\le 1\}$, where
$\tau = s b$ and $g_{a+\tau }$ is the symplectic reduction of $g$ on $X_{a+\tau}$.

We can now prove:
\begin{theorem}\label{theorem:twistedricciflat}
Let $(M,g)$, $h(\tau)$ etc. be as above. Suppose that for some function $f = f(\tau)$ on $M$,
$g$ satisfies the following equation:
\begin{equation}\label{twistedricciflat}
{\rm Ric(g)} +\frac{\sqrt{-1}}{2} \p \bar \p \left (\log \det (w^{ij})+ f\right )=\lambda \omega_g.
\end{equation}
Then we have:

\noindent
(1) The function
$\Delta _g \tau _k -\frac {\p w^{kl}}{\p \tau _l} - w^{kl} \frac {\p f}{\p \tau _l}$
is constant along each connected component of $\mu^{-1}(a+\tau)$;

\noindent
(2) Either ${\rm Ric}(h(0)) = \lambda \omega_{h(0)}$, i.e., $h(0)$ is K\"ahler-Einstein, or
$h(\tau)=\phi _{\tau } ^* g_{a+\tau}$ is a solution of the normalized K\"ahler-Ricci
flow on $X_a$:
\begin{equation}
\label{eq:krflow}
\frac{\p \omega}{\p t} = - {\rm Ric}(\omega) + \lambda \omega,
\end{equation}
provided that $\tau_k(t) = {c_k} \left ( e^{\lambda t} -1\right )/ {\lambda}$ {\rm ($1\le k\le N$)}
\footnote{If $ \lambda=0$, then $\tau_k(t)= c_k t $}, where
\begin{equation}\label{time}
c_k = \left( -\frac{1}{4} \Delta _g \tau _k +
\frac{1}{4} \frac {\p w^{kl}}{\p \tau _l}  + \frac{1}{4} w^{kl} \frac
{\p f}{\p \tau _l}\right)\big |_{\mu^{-1}(a)}.
\end{equation}

\end{theorem}
\begin{proof}
Since $\log \det (g)= \log \det (h) + \log \det (w^{ij})$ in a certain local holomorphic frame, we see that
$${\rm Ric }(g) + \frac{\sqrt{-1}}{2} \p \bar \p \left(\log \det (w^{ij}) + f \right)=\lambda \; \omega_g$$
is equivalent to:
$$\frac{\sqrt{-1}}{2} \p \bar \p \left( -\log \det (h)  +f\right) =\lambda \; \omega_g.$$
In turn, by Lemma \ref{complexhessian} and the assumption that $f=f(\tau)$, we show that the above equation is equivalent to the
following system:

\begin{equation}\label{equivalentsystem}
\left\{ \begin{aligned} &(i) \;\;\; 4 Ric (h)_{i\bar j} - w^{lk}
    \left( \frac{\p \log \det (h)}{\p \tau _l}-  \frac{\p  f}{\p \tau
      _l}\right) \frac{\p h _{i\bar j}}{\p \tau _k} = 4 \lambda \;  h _{i\bar
    j}\\ & (ii) \;\; \; \frac{\p}{\p z_k}\left(w^{lk}
    \left( \frac{\p \log \det (h)}{\p \tau _l}- \frac{\p  f }{\p \tau
      _l}\right)\right) =0\\& (iii)\;\;\; \frac{\p}{\p \bar z_k}\left(w^{lk}
    \left( \frac{\p \log \det (h)}{\p \tau _l}-  \frac{\p  f }{\p \tau
      _l}\right)\right) =0 \\
&(iv) \;\; \;  \frac{\p}{\p \tau _k}\left( w^{lj}
    \left( \frac{\p \log \det (h)}{\p \tau _j}-  \frac{\p  f }{\p \tau
      _j}\right)\right) = - 4 \lambda \delta _{kl}
\end{aligned}\right.
\end{equation}
\noindent
It follows form (ii) and (iii) in the system \eqref{equivalentsystem} that
$$ w^{lk} \left( \frac{\p \log \det (h)}{\p \tau _l}-  \frac{\p  f}{\p \tau _l}\right)$$
is a function of the Hamiltonian coordinates $\tau_1, \cdots, \tau_N$ only.

On the other hand, from Lemma \ref{lemma:hamiltonianlaplacian}, one can infer that:
$$\frac{1}{4}  w^{kl}  \left( \frac{\p \log \det (h)}{\p \tau _l}-  \frac{\p  f}{\p \tau
_l}\right) =\frac{1}{4} \Delta _g \tau _k -
\frac{1}{4} \frac {\p w^{kl}}{\p \tau _l}  -\frac{1}{4} w^{kl}
\frac{\p  f}{\p \tau  _l}.$$
This shows (1).

It follows from (iv) in \eqref{equivalentsystem}, Lemma \ref{lemma:hamiltonianlaplacian} and \eqref{time} that
\begin{equation}
\label{eq:derivativeof tau}
\frac{1}{4} w^{lk} \left( \frac{\p \log \det (h)}{\p \tau _l}-  \frac{\p  f}{\p \tau _l}\right) = -\lambda \tau_k - c_k.
\end{equation}
By our choice of $\tau_k(t)$, this implies
\begin{equation}
\label{timechange}\frac{\p \tau _k }{\p t}= -\frac{1}{4} w^{lk}
\left( \frac{\p \log \det (h)}{\p \tau _l}-  \frac{\p  f}{\p \tau
_l}\right).\end{equation}
Observe that $\tau_k'(t)\not= 0$ for all $t$ whenever $c_k\not= 0$.
Hence, if $h(0)$ is not a K\"ahler-Einstein metric, then $\tau(t)$ is a genuine parameter change of time $t$.
Thus we have derived from (i), (ii) and (iii) of (\ref{equivalentsystem}) the K\"ahler-Ricci flow
$$\frac{\p h}{\p t} = -{\rm Ric}(h) +\lambda\,h,~~~~{\rm on}~X_a.$$

\noindent
\end{proof}

\noindent
Let $(M,g)$, $h$ be in the above theorem. If $h(0)$ is not K\"ahler-Einstein, then it follows from \eqref{timechange} and a direct computation
\begin{equation}\label{positivity}
\begin{aligned} &\, R-n\lambda +\frac{\p f}{\p t} \\
= &\, \frac{\p f}{\p t} -\frac{\p \log \det (h)}{\p t}\\
 =&\, \sum _k \left( \frac{\p f}{\p \tau _k}- \frac{\p \log \det (h)}{\p
\tau _k} \right) \frac{d\tau _k } {d t}\\
=&\, \frac{1}{4}\sum _{k,l} w^{kl} \left(\frac{\p f}{\p \tau _k}-  \frac{\p \log \det (h)}{\p
\tau _k} \right)\left( \frac{\p f}{\p \tau _k}- \frac{\p \log \det (h)}{\p
\tau _l}\right) .\end{aligned}\end{equation}
Since $w_{ij}$
(and hence $w^{ij}$) is positive definite, i.e., $w^{kl} \xi _l \xi
_k>0$ for every non-zero $(\xi _1, \cdots, \xi _N)$, we have
\begin{equation}
\label{eq:df-ineq}
R(h) - \lambda n + \frac{\p f}{\p t} \ge 0
\end{equation}
and the equality holds at some $t$ if and only if $h$ is K\"ahler-Einstein, where $R(h)$ denotes the scalar curvature of $h$.

There is an integral condition on the descended solution $h(t)$ of the K\"ahler-Ricci flow from a solution of \eqref{twistedricciflat}:
For simplicity, we assume that $N=1$, that is, the action group is $S^1$.
By \eqref{closed}, we have
$$d \theta |_{\mu^{-1}(a+\tau)}  = - \frac{\sqrt{-1}}{2} \frac{\partial h(\tau)}{\partial \tau }=
\frac{\sqrt{-1}}{2} \frac{dt}{d\tau} \left ({\rm Ric}(h(\tau))- \lambda \omega_{h(\tau)}\right).$$
On the other hand, using the K\"ahler-Ricci flow, we can show
$$-c_1(X_a) + \lambda [\omega_{h(t)}] = e^{\lambda t} \left (-c_1(X_a) + \lambda [\omega_{h(0)}]\right ).$$
Here $[\omega]$ denotes the cohomology class represented by $\omega$. It follows from the above equations
$$ d \theta |_{\mu^{-1}(a+\tau)} = \frac{1}{2 c} \left (c_1(X_a) - \lambda [\omega_{h(0)}]\right ).$$
Noticing that $\theta|_{\mu^{-1}(a+\tau)}$ is a connection of the circle bundle $\pi: \mu^{-1}(a+\tau) \mapsto X_{a+\tau}$, so its curvature
$d\theta$ represents its first Chern class. Hence, $\lambda [\omega_{h(0)}]$ must be in $H^2(X_a, 2\pi \ZZ)$.
In particular, if $\lambda=0$, the above shows that
the associated circle bundle is just the pluri-anti-canonical bundle.

The converse of the above theorem is given in the following:
\begin{theorem}\label{theorem:converse}
Let $X$ be a K\"ahler manifold. If $\tilde h$ is a solution of (\ref{eq:krflow}) on $ X \times [t_0,t_1]$
such that $\lambda [\omega_{\tilde h(0)}]$
lies in $H^2(X, 2\pi \ZZ)$, then there is
a unique principal $S^1$-bundle over $ X \times [t_0,t_1]$ and a $S^1$-invariant metric $g$ on $M$ satisfying
the equation (\ref{twistedricciflat}) for some $f$ and a function $\tau(t): [t_0,t_1]\to \RR$
such that $\tau(t_0)=0$ and $\tilde h(t)= h(\tau(t))$, where $h(\tau)$ is induced from $g$ as in last theorem. Moreover, the
curvature of the principle bundle $M$ is given by
$$\gamma _k :=\sqrt{-1}\left\{ -\frac{1}{2} \frac{\partial {h_{i\bar
    j}}}{\partial \tau } dz_i\wedge d \bar z _j - \frac {\partial
  w }{\partial z_j} d\tau \wedge dz_j + \frac {\partial
  w }{\partial \bar z_j}d\tau \wedge d\bar z_j\right \}.$$
\end{theorem}
\begin{proof}
First we assume that $\tilde h$ is K\"ahler-Einstein, i.e., ${\rm Ric}(\tilde h) = \lambda \tilde h$. Take $M=X\times \CC$
and the vector field $V= 2\,{\rm Im}(z\frac{\partial}{\partial z})$ is simply
the one inducing the standard rotation on $\CC$. The lifting metric $g$ is of the form
$$g = \tilde h + w d\tau^2 + w^{-1} \theta^2,$$
where $\theta $ is the dual of $V$ as we defined before. Then $\frac{1}{2} |z|^2$ is the associated moment map. Define
$w^{-1} = |z|^2= 2  \tau^2 $ and $f$ as a function of $\tau$
by
$$\frac{\p f}{\p \tau} \,=\, 4\lambda\tau\,w.$$
Then one can check directly that $g$ satisfies
$${\rm Ric}(g) + \frac{\sqrt{-1}}{2} \partial \overline{\partial} \left (\log |V|^2 + f\right ) = \lambda \omega_g.$$
Hence, we get a lifting of $\tilde h$ on $X\times S^1\times [\tau_0, \tau_1]$, where $\tau_0$ and $\tau_1$ are determined by $t_0$ and $t_1$,
respectively.

Now we suppose that $\tilde h$ is a non-static solution of the K\"ahler-Ricci flow:
$$\frac{\p \tilde h}{\p t} = -{\rm Ric }_{X}(\tilde h) + \lambda \,\tilde h.$$
We want to find $(M,J)$ and a K\"ahler metric $g$ of the form
$$\tilde h_{i\bar j} dz_i dz_j + w d\tau^2 + w^{-1} \theta ^2 ,$$
which satisfies \eqref{twistedricciflat}:
$${\rm Ric }(g) +\frac{\sqrt{-1}}{2} \p\bar \p (- \log w + f) = \lambda \omega_g$$
for some $f=f(t)$, i.e., it is constant along $X$. We will revert the reduction of $g$ to $h(\tau)$ in the previous theorem.

As expected, we set
$$\tau (t) = \frac{c}{\lambda} (e^{\lambda t}-1),$$
where $c$ is a positive constant. Hence,
$$\frac{d \tau}{dt} = c e^{\lambda t}.$$

Choose any smooth function $f=f(\tau)$ such that
$$ -R(\tilde h(t)) + n\lambda -  \frac{\p  f(\tau(t))}{\p t} < 0$$
on $X\times [t_0,t_1]$, where $R(\tilde h(\tau))$ is the scalar curvature of $\tilde h(\tau)$. Now we define
$$ w = \frac{  e^{- 2\lambda t} \left (R(\tilde h(t)) - n\lambda + \frac{\p  f(\tau(t))}{\p t} \right )}{4  c^2} >0.$$
Since $\tilde h$ is a solution of the K\"ahler-Ricci flow, we have
$$ \frac{\p \log \det (\tilde h)}{\p t} = -R(\tilde h(t)) + n\lambda.$$
If we regard $\tilde h$ as a function of $\tau$, we can deduce from the above that:
\begin{equation}\label{timew}
\frac{\p \tau  }{\p t} = -\frac{1}{4} w^{-1}
\left( \frac{\p \log \det (\tilde h)}{\p \tau }-  \frac{\p  f}{\p \tau}\right).
\end{equation}
\noindent
Define a 2-form as follows:
\begin{equation}\label{dtheta}\gamma  :=\sqrt{-1}\left\{ -\frac{1}{2} \frac{\partial {\tilde h_{i\bar
    j}}}{\partial \tau } dz_i\wedge d \bar z _j - \frac {\partial
  w }{\partial z_j} d\tau \wedge dz_j + \frac {\partial
  w}{\partial \bar z_j}d\tau \wedge d\bar z_j\right \}.
  \end{equation}
We claim that $\gamma$ is closed on $X\times [t_0,t_1]$. The closedness of $\gamma$ is equivalent to
\begin{equation}\label{closedness}
 \frac {\p ^2 \tilde h_{k\bar l}} {\p \tau^2} = - 4\frac{\p ^2 w}{\p z_k \p \bar z_l}.
 \end{equation}
Using the Ricci flow and the definition of $\tau(t)$, we see that the left-handed side becomes
$$\frac {\p}{\p \tau}\left (c^{-1} e^{-\lambda t}\left (- {\rm Ric}( \tilde h)_{k\bar l} + \lambda \tilde h_{k\bar l}\right )\right )
= - c^{-1} e^{-\lambda t} \frac{\p^2 }{\p z_k\p \bar z_l}\left (\frac{\p \log \det (\tilde h)}{\p \tau}\right ).$$
On the other hand, using \eqref{timew}, we have
$$ 4 w = - c^{-1} e^{-\lambda t} \left (\frac{\p \log \det (\tilde h)}{\p \tau}- \frac{\p f}{\p \tau}\right ).$$
The claim follows.

Then we can take $M$ to be the unique principal $S^1$-bundle  $\pi: M \to X \times [t_0,t_1]$ with the connection
1-form $\theta$ such that
$$d\theta  = \pi ^* \gamma .$$
Here we have used the assumption on the K\"ahler class of $\tilde h(0)$.

We then define the complex structure $J$ on $M$ by imposing that $J \theta  = - w d\tau $ and it restricts to the given one on $X$.
Since $T^{1,0}M$ is locally spanned by $dz_1,\cdots, dz_n$ and $w d\tau -\sqrt{-1} \theta$, $J$ is integrable if
$d (w d \tau - \sqrt{-1} \theta)$ has no (0,2)-components. The latter can be checked directly by using the definition of $\theta$.
Hence, $J$ is integrable. Furthermore, we can endow $(M,J)$ with a K\"ahler structure:
By the definition, we have
$$d \theta  = \sqrt{-1}\left\{ -\frac{1}{2} \frac{\partial {\tilde h_{i\bar
    j}}}{\partial \tau } dz_i\wedge d \bar z _j - \frac {\partial
  w }{\partial z_j} d\tau \wedge dz_j + \frac {\partial
  w }{\partial \bar z_j}d\tau\wedge d\bar z_j\right \}.$$
It follows that $\omega _{\tilde h} - d\tau \wedge \theta $ is a closed 2-form.
Clearly, this is the K\"ahler form of the required K\"ahler metric
$$g= \tilde h_{i\bar j} dz_i dz_j + w d\tau^2 + w^{-1} \theta ^2,$$
that is, $\omega_g = \omega _{\tilde h} - d\tau  \wedge \theta $.

Let $V$ be the vector field inducing the standard clock-wise rotation on the circle bundle $M$, then
$\theta (V)=1$ and $i_V\omega_g = d\tau$. This means that $\tau$ is a moment map.
from the construction, we can easily show that $\tilde h$ coincides with $h(\tau)$ from last Theorem.

\end{proof}

\begin{remark}
The above lifting is not unique since we do have choices of $f$. If $g$ and $g'$ are such metrics corresponding to
$f$ and $f'$, respectively, then we notice that $\gamma$ is independent of $f$ and $f'$, so we have the same circle bundle
$M$. Moreover, the symplectic form $\omega_g$ is independent of the choice of $f$.

One may replace \eqref{twistedricciflat} by a slightly more general equation: From
the above proof, one can see that any solution of
${\rm Ric}(g)+ \p\bar \p (\log \det (w^{ij})+f)= \Omega$ also descends to a solution of the
K\"ahler-Ricci flow, so long as $\Omega$ is a closed $(1,1)$-form such
that $(\pi _{a+\tau}) _* \Omega = \lambda\;
h(\tau)$ and $\Omega (Z,V_i )=0$ for
every $Z\in Q(a+\tau)$ for every $\tau$. Of course, it holds for $\Omega =\lambda \;\omega _g$.
\end{remark}

In the case in which the action group is just $S^1$, \eqref{time}
takes a particularly interesting form as it reduces to a
(modified) mean curvature flow:
\begin{lemma}\label{meancurvature}
One has that
$$w^{-1}\,\frac{\p \log \det (h)}{\p
  \tau}= - \frac{H(\tau)}{|V|_g}  - \frac{1}{2} \frac{\p w^{-1} }{\p
  \tau}  $$
where $H(\tau)$ is the mean curvature of $\mu^{-1}(a+\tau)$ with respect to the unit normal $\sqrt{w} JV$.
\end{lemma}
\begin{proof}
Given any smooth function $f$ on a Riemannian manifold $(N,g)$, along its level set $N_a:=\{
x\in M : \; f(x)=a\}$, we have
$$ Hess _f(Y_1,Y_2) = -\langle (\nabla _{Y_1} Y_2) , \nabla f\rangle =
-\langle (B(Y_1, Y_2) , \nabla f\rangle,$$
where $Y_1,Y_2$ are tangent to $N_a$ and $B$ denotes the 2nd fundamental form of $N_a$.
It follows
\begin{equation}\label{hessinatomeancurvature}
\Delta _g f \mid _{N_a} = -\langle H ,\nabla f\rangle_g + Hess _f (\nu, \nu),\end{equation}
where $H$ is the mean curvature of $N_a$ and $\nu=\frac{\nabla f}{|\nabla f|}$ is the unit normal.

Now applying \eqref{hessinatomeancurvature} to the moment map $\tau$ regarded as a function on $M$, we get
$$\Delta _g \tau = -\langle H, \nabla \tau \rangle _g  + Hess
_{\tau} (\nu , \nu).$$
On the other hand, since  $\nabla \tau = JV$ and $\nu = \sqrt{w} \nabla \tau$, a straightforward calculation
shows
$$Hess_\tau(\nu ,\nu)= \nu(\nu \tau) = \nu (\sqrt{w} d\tau (JV)) = \nu ( w^{-1/2}) = w^{-1/2} \frac{\p w^{-1/2}}{\p \tau}
= \frac{1}{2} \frac{\p w^{-1} }{\p
  \tau} .$$
Then the claim follows from \eqref{hamiltonianlaplacian}.

\end{proof}

It follows from this lemma that the derivative $\frac{ d\tau}{dt}$ in Theorem \ref{theorem:twistedricciflat} satisfies
an evolution equation of mean curvature flow type:
$$
\frac{ d\tau}{dt} = \frac{H(\tau)}{ 4 |V|_g}  + \frac{1}{8} \frac{\p w^{-1} }{\p \tau} + \frac{1}{4} w^{-1} \frac{\p f}{\p \tau}.$$

\section{Scalar V-soliton equation}\label{continuity}

In this section, we will address the solvability of the following complex Monge-Ampere equation:
\begin{equation}\label{twistedmongeampere}
(\omega_{g_0} +\frac{\sqrt{-1}}{2} \p\bar \p  u )^n = \left(|V| _{g_0}^2 + \frac{\sqrt{-1}}{2} \p\bar \p u (V,JV) \right)
e^{F-\lambda u}\,\omega_{g_0}^n,
\end{equation}
where $g_0$ is a given K\"ahler metric and $F$ is a given function satisfying
$$\int_M \left(|V|^2_{g_0} e^F - 1\right) \omega_{g_0}^n = 0. $$
We call \eqref{twistedmongeampere} scalar V-soliton equation.
We will assume that both $g_0$ and $F$ are invariant under the
$S^1$-action induced by $V$.

Our main goal here is to develop some preliminary estimates
necessary to prove the existence of solutions for this scalar V-soliton equation.
Higher order estimates will be done in a forthcoming paper.
For simplicity, we assume
that $M$ is compact.

One motivation for studying (\ref{twistedmongeampere}) comes from establishing the existence of V-solitons:
Suppose that $N=1$ and $g$ is a solution of (\ref{twistedricciflat}), that is, $g$ is a V-soliton metric.
Now we choose $g_0$ such that $c_1(M)$ coincides $\lambda [\omega_{g_0}]$. then we can write $g$ as
$$\omega_g = \omega_{g_0} + \frac{\sqrt{-1}}{2} \p \bar\p u.$$
Define $F=F'+f$ and $F'$ by the equation:
$${\rm Ric}(g_0)-\lambda \omega_{g_0} = \frac{\sqrt{-1}}{2} \p \bar \p F'.$$
Such a function $F'$ is only determined up to constants. Then $g$ is a $V$-soliton metric if and only if $u$ satisfies
(\ref{twistedmongeampere}) modulo addition of constants.
In fact, one can easily see that even if $M$ is not compact, if $u$
is a solution of (\ref{twistedmongeampere}), then $g$ defined as above in terms of $u$
is still a $V$-soliton.

Here we consider \eqref{twistedmongeampere} only when $\lambda =0$. The case for $\lambda=-1$ can be done in a similar and simpler way.
As usual, the case for positive $\lambda$ is more tricky.
We will assume that $u$ is invariant.

\begin{lemma}\label{boundednessofV(u)}
There is a uniform constant $C=C(g_0)$ such that
for any $S^1$-invariant function $u$ with $\omega_{g_0} + \frac{\sqrt{-1}}{2} \p \bar\p u\ge 0$,
we have
$$ w^{-1}_0 \left |\frac{\p u}{\p \tau}\right | \le C,$$
where $w_0^{-1} = |V|_{g_0}^2$.
\end{lemma}
\begin{proof} This is a known fact (cf. \cite{zhu}). For the readers' convenience, we include a sketched proof here.
As before, we denote by $\mu=\tau$ the moment map associated to the $S^1$-action by $V$. Since M is compact, $\mu$ has at least
two critical points, so $V$ has at least two zeroes. To estimate $JV(u)=w^{-1}_0 \left |\frac{\p u}{\p \tau}\right |$ at any given $p\in M$, we pick up
a trajectory $\gamma$ of the gradient $\nabla \mu$ from one critical point to another. Since $w^{-1}_0=0$ at critical points of $\mu$, we may assume
$p$ is not a critical point. Then $\gamma$ sweeps out a holomorphic sphere $S$ with two punctures by the $S^1$-action.
Those two punctures are exactly those critical points which $\gamma$ connects. Using the $S^1$-symmetry, we get
$$\omega_{g_0} (V, JV) + \frac{\sqrt{-1}}{2} \p \bar\p u (V,JV) > 0   ~~~~~~~~~~~ {\rm on} ~ S.$$
In view of \eqref{hessian}, this is the same as
$$\frac{\p }{\p \tau} \left (JV (u)\right ) = \frac{\p }{\p \tau} \left (w_0^{-1} \frac{\p u}{\p \tau}\right ) > - 4.$$
Integrating this along $\gamma$ starting from either $\tau_{max}=\sup_\Gamma \mu$ or $\tau_{min}=\inf_\gamma \mu$, we get
$$-4(\mu(p)-\tau_{min}) \le JV(u)(p) \le 4(\tau_{max} - \mu(p)).$$
It follows that $|JV(u)|\le 4(\tau_{max}-\tau_{min})$, so the lemma is proved.
\end{proof}
We may use the perturbation method to solve \eqref{twistedmongeampere}. Consider
\begin{equation}\label{perturbation}
(\omega_{g_0} +\frac{\sqrt{-1}}{2} \p\bar \p  u )^n = \left(\epsilon + |V| _{g}^2 \right)
e^{F+c_\epsilon}\,\omega_{g_0}^n,
\end{equation}
where $\epsilon > 0$, $\omega_g =\omega_{g_0} + \frac{\sqrt{-1}}{2} \p\bar \p  u$ and $c_\epsilon$ is chosen such that
$$\int_M \left((\epsilon + |V|^2_{g_0}) e^{F_\epsilon} - 1\right) \omega_{g_0}^n =0,$$
where $F_\epsilon = F + c_\epsilon$.

Now let us introduce some notations. Set
$$C^{k,\alpha} (M,V) :=\left\{u\in C^{k,\alpha } (M)\,| \; \; V(u) = 0 \right\},$$
where $C^{k, \alpha } (M)$ is the H\"older space of $C^k$-smooth functions such that
$$|| u ||_{ C^{k, \alpha  }}:=  \sum _{i=1}^k \sup _{x\in M}
|\nabla ^i u|+ \sup _{x, y\in M, \; x\neq y} \frac{ |\nabla ^k u (x) - \nabla ^k u (y)|}{d(x,y)^{\alpha}}<+\infty,$$
where $d(\cdot, \cdot)$ denotes the distance function of any fixed metric $g$.
Clearly, this coincides with the space $C^{k, \alpha} (M)_{S^1}$ which consists of $S^1$-invariant functions in $C^{k,\alpha }(M)$.
We further set
$$C^{k, \alpha } (M;V)_g :=\left\{v\in C^{k, \alpha } (M;V)\, | \; \
  \int _M \left((\epsilon+|V|_{g}^2) e^v - 1\right) \omega_{g}^n=0\right\}$$
and
$$C^{k, \alpha  }_g (M;V) :=\left\{u\in C^{k, \alpha } (M;V)\, |\; \
  \int _M u \,  \omega
_{g} ^n=0\right\}.$$
For $k\ge 2$, we also denote by $P^{k,\alpha}(M,V)$ the set of all $u\in C^{k,\alpha}(M,V)$ such that
$\omega_u:= \omega_{g_0} +\frac{\sqrt{-1}}{2} \p \bar\p u > 0$.
Define a differential operator from $P^{k,\alpha}(M,V)$:
$$\Phi_\epsilon (u) := \log \left(\frac{(\omega_{g_0} +\frac{\sqrt{-1}}{2} \p \bar\p u)^n}{\omega_{g_0}^n}\right) -\log (\epsilon + |V|_u^2), $$
where $|V|_u^2 = \omega_u (V,JV)$ is the square norm of $V$ with respect to the metric given by $\omega_u$.

Clearly, for $k\ge 2$, $\Phi_\epsilon$ maps into $C^{k-2, \alpha } (M;V)_{g_0} $.
To solve \eqref{perturbation}, we only need to show that $\Phi_\epsilon$ is surjective.
We will prove that for $k$ sufficiently large,\footnote{$k\ge 4$ should be sufficient.}
$$\Phi_\epsilon (P^{k,\alpha}(M,V)\cap C^{k, \alpha  }_{g_0} (M;V)) = C^{k-2, \alpha } (M;V)_{g_0}.$$

The tangent space to $C^{k-2, \alpha } (M;V)_{g_0}$ at $\Phi+\epsilon(u)$ is the space:
$C^{k-2, \alpha  }_{g_u}(M;V )$. Hence, the differential $D\Phi _\epsilon|_{u}$ of $\Phi_\epsilon$ at $u$ is a
linear map from $C^{k, \alpha } (M;V)$ into $C^{k-2,\alpha}_{g_u}(M;V)$. Furthermore, we have

\begin{lemma}\label{differential}
For any $\epsilon > 0$, $\Phi_\epsilon$ is an elliptic operator. Moreover, for any $u\in P^{k,\alpha}(M,V)$,
the differential $D\Phi _\epsilon|_{u}$ is surjective with only constant functions in its kernel.
\end{lemma}
\begin{proof} The ellipticity of $\Phi_\epsilon$ means that
for any $u\in P^{k,\alpha}(M,V)$, $D\Phi _\epsilon|_{u}$ is elliptic. A straightforward computation shows:
$$D\Phi _\epsilon|_{u}(\bar{u})=\Delta _{g_u} \dot{u} -\frac{\sqrt{-1}}{2(\epsilon+|V|_{u}^2)}
\,\p \bar\p \dot {u} (V,JV).$$
At any given point $p\in M$, we can choose a basis $\{e_i\}$ of $T^{1,0}_pM$ satisfying:
$$ g_u(e_i, \bar e_j) = \delta_{ij},~~~~\frac{\sqrt{-1}}{2} \p \bar\p u ( e_i, \bar e_j) = a_i \delta_{ij}.$$
In terms of this basis, we have
$$D\Phi _\epsilon|_{u}(\bar{u})(p)= \sum_i \frac{(\epsilon + \sum_{j\not= i} |v_j|^2) a_i }{\epsilon + \sum_j |v_j|^2}.$$
This shows the ellipticity of $D\Phi _\epsilon|_{u}$ at $p$, and consequently, ellipticity of $\Phi_\epsilon$.

Moreover, it follows from the above computation and the Maximum principle that $D\Phi _\epsilon|_{u}$ is surjective
and its kernel consists of only constant functions.

\end{proof}

\begin{remark}
In fact, one can also show that $D\Phi _\epsilon|_{u}$ is self-adjoint. Moreover, by the above, we see that
$D\Phi _0|_{u}$ is also elliptic but degenerate.
\end{remark}

As said, we may use the continuity method to solve \eqref{perturbation}.
Fix a large $k>0$. Choose any path $F_{\epsilon,s}$ in $C^{k-2, \alpha } (M;V)_{g_0}$ ($s\in [0,1]$)
with $F_{\epsilon,0}=-\log (\epsilon + |V|^2_{g_0})$ and $F_{\epsilon,1}$ coincides with $F_\epsilon$
in \eqref{perturbation}. Consider a family of complex Monge-Ampere equations:
\begin{equation}\label{continuity}
(\omega_{g_0} +\frac{\sqrt{-1}}{2} \p\bar \p  u )^n = \left(\epsilon + |V| _{g}^2 \right)
e^{F_{\epsilon,s}}\,\omega_{g_0}^n,
\end{equation}
where $\omega_g =\omega_{g_0} + \frac{\sqrt{-1}}{2} \p\bar \p  u$. Define
$$I=\{s\in [0,1]~|~ \eqref{continuity}~{\rm has~a~solution~for ~any}~s'\in [0,s]\}.$$
Clearly, $0\in I$ since $u=0$ is a solution. It follows from the above lemma and the Inverse Function Theorem
\begin{corollary}\label{open}
The set $I$ defined above is open.
\end{corollary}

Hence, establishing the existence of a solution for \eqref{perturbation} is equivalent to proving that $I$ is closed.
For this purpose, we need a prior estimates for solutions of \eqref{continuity}.
In view of the $C^{2,\alpha}$-estimate due to Evans, Krylov etc.. (cf. \cite{ev}, \cite{kr}, \cite{ca}, \cite{tr},
 \cite{ti1}), it suffices to have a priori $C^2$-estimates for \eqref{continuity}.

\subsection{The $C^0$-estimate}

The purpose of this subsection is to derive a $C^0$-estimate by the standard Moser iteration,
keeping track of the dependence on $\epsilon> 0$. First, we have
\begin{proposition}\label{c0}
There is a uniform constant $C$ which depends only on $(M,g_0)$, $||F||_{C^1(M)}$,
$\sup_M |V|_{g_0}$ and $\sup_M |{\rm div}(JV)|$
such that for any solution $u$ of \eqref{continuity} with $\int_M u\, \omega_{g_0}^n =0$, we have
$$\sup _M |u| \leq C. $$
\end{proposition}
\begin{proof} First we assume that $u$ is a solution of \eqref{continuity} for some $s\in [0,1]$ such that $\sup_M u =-1$.
For simplicity, denote $F_{\epsilon,s}$ by $\bar F$ and $u_-=- u$. Note that $u_-\ge 1$.
Integrating by parts, we get for $\ell \ge 1$
\begin{equation}\label{intbyparts}
\begin{aligned}
&~~n\,\int_M \; u_-^\ell \,( \omega_{g} ^n  -\omega _{g_0}^n )\\
& =  \frac{n \ell \sqrt{-1} }{2} \int_M \;
u_-^{\ell-1} \,\left( \p u_- \wedge \bar\p u_-\wedge \sum_{j=0}^{n-1} \omega_{g_0}^j \wedge \omega _{g}^{n-j-1}\right)\\
& \geq \ell \int _M
u_-^{\ell-1} |\nabla _{g_0} u_-| ^2\, \omega_{g_0}^n \\
&= \frac{4\ell}{(\ell +1)^2}\,\int _M |\nabla _{g_0}
u_-^{\frac{\ell+1}{2}}|^2 \,\omega^n_{g_0}\end{aligned}
\end{equation}
Multiplying $u_-^\ell$ on both sides of \eqref{continuity} and integrating, we deduce from the above
\begin{equation}
\begin{aligned}
\int _M |\nabla _{g_0} u_-^{\frac{\ell +1}{2}}|^2 \omega^n_{g_0} &\leq \frac{n
(\ell + 1) ^2}{4 \ell}
\int_M u_-^\ell\, \left( \omega _{g}^n - \omega
  _{g_0} ^n\right) \\&  = \frac{n (\ell + 1)^2}{4 \ell}
 \int_M u_-^\ell \,\left( (\epsilon + |V|_g^2)\, e^{\bar F} - 1 \right)  \omega_{g_0}^n
 \end{aligned}
\end{equation}
Using Lemma \ref{complexhessian} and noticing $w_0\,JV = \frac{\p}{\p \tau}$, we can compute
$$\frac{\sqrt{-1}}{2} \p\bar\p u (V,JV) = \frac{1}{4} \,w_0^{-1}\,\frac{\p}{\p \tau } \left( w_0^{-1}\,\frac{\p \phi}{\p \tau }
\right)= \frac{1}{4} \, JV(JV(u)).$$
It follows
\begin{equation}
\label{c0-estimate-1}
\int _M |\nabla _{g_0} u_-^{\frac{\ell +1}{2}}|^2 \omega^n_{g_0}\,\le \,n \ell
\int_M u_-^\ell\,\left( e^{\bar F} (\epsilon + |V|_{g_0}^2 + \frac{1}{4} JV(JV(u))) - 1 \right )   \omega _{g_0}^n.
\end{equation}
Write $W=JV- \sqrt{-1} V$. Then $W$ is a holomorphic vector field. Since $V(u)=0$, we have $W(u)=JV(u)$ is real-valued and bounded.
Recall the identity:
$${\rm div}(u_-^\ell e^{\bar F} W(u) W) = {\rm div}(W)u_-^\ell e^{\bar F} W(u) + W(u_-\ell e^{\bar F} W(u)),$$
where the divergence ${\rm div}(W)$ is taken with respect to the metric $g_0$.
Therefore, there is a constant $C_1$ which depends only on $(M,g_0)$, $||F||_{C^1(M)}$,
$\sup_M |V|_{g_0}$ and $\sup_M |{\rm div}(JV)|$ such that
$$u_-^\ell e^{\bar F} W(W(u))\, \le \,{\rm div}(u_-^\ell e^{\bar F} W(u) W) + \ell u_-^{\ell-1} |W(u)|^2 e^{\bar F}
+ C_1 \, u_-^\ell |W(u)|.$$

Plugging this into \eqref{c0-estimate-1} and using Lemma \ref{boundednessofV(u)}, we obtain
\begin{equation}
\label{c0-estimate-2}
\int _M |\nabla _{g_0} u_-^{\frac{\ell +1}{2}}|^2 \omega^n_{g_0}\,\le\, C_2 \,\ell \,\int_M u_-^{\ell-1}(u_- + \ell ) \omega_0^n,
\end{equation}
where $C_2$ is a uniform constant depending only on $(M,g_0)$, $||F||_{C^1(M)}$,
$\sup_M |V|_{g_0}$ and $\sup_M |{\rm div}(JV)|$.

Since $u_-\ge 1$, it follows
\begin{equation}
\label{c0-estimate-3}
\int _M |\nabla _{g_0} u_-^{\frac{\ell +1}{2}}|^2 \omega^n_{g_0}\,\le\, 2\,C_2 \,\ell^2 \,\int_M u_-^{\ell +1} \omega_0^n.
\end{equation}

Now we can apply the standard
Moser iteration scheme: Denote by $C_S$ the Sobolev constant for $g_0$, then for any smooth function $f$ on $M$, we have
(cf. \cite{gt} Theorem 7.10):
$$|| f ||_{\frac{2n}{n-1}} \,\leq\, C_S\, \left( || \nabla f||_2 +
||f||_2\right).$$
Applying this to \eqref{c0-estimate-3} for $f=u_-^{\frac{\ell+1}{2}}$, we get
\begin{equation}\label{c0-estimate-4}
\left( \int_M  u_-^{\frac{ n (\ell +1)}{n-1}}\,\omega_{g_0}^n\right)
^{\frac{n-1}{n(\ell +1)}}\, \leq\, \left (C_3\, \ell^2 \, \int _M u_-^{\ell+1} \omega_{g_0}^n\right )^{\frac{1}{\ell +1}},
\end{equation}
where $C_3$ is a uniform constant depending only on $(M,g_0)$, $C_S$, $||F||_{C^1(M)}$,
$\sup_M |V|_{g_0}$ and $\sup_M |{\rm div}(JV)|$.

Set $\ell_1=1$ and $\ell_{i+1} =\frac{2n}{n-1} (\ell_i+1) -1$ inductively for $i\ge 1$. Then
we have
$$\left( \int_M  u_-^{\ell_i +1}\,\omega_{g_0}^n\right)
^{\frac{1}{\ell_i +1}}\, \leq\, \prod _{j=1}^{i-1} \left (C_3\,( \ell_{j}+1)^2 \right )^{\frac{1}{\ell_j+1}}
\, \left(\int _M u_-^2 \omega_{g_0}^n\right )^{\frac{1}{2}}.$$
Note that
$\ell_{j}+1= 2 \left (\frac{n-1}{2n}\right)^j$, we can deduce from the above
\begin{equation}
\label{c0-estimate-5}
\sup_M u_- \,\le\, C_4 \left(\int_M u_-^2 \,\omega_{g_0}^n\right )^{\frac{1}{2}},
\end{equation}
where $C_4$ is a uniform constant depending only on $(M,g_0)$, $C_S$, $||F||_{C^1(M)}$,
$\sup_M |V|_{g_0}$ and $\sup_M |{\rm div}(JV)|$.

Moreover, applying the Poincare inequality to \eqref{c0-estimate-2} with $\ell=1$ and noticing $u_-\ge 1$, we get
\begin{equation}
\label{c0-estimate-6}
\left(\int_M u_-^2 \,\omega_{g_0}^n\right )^{\frac{1}{2}}\, \le\, C_5 \int_M u_- \, \omega_{g_0}^n,
\end{equation}
where $C_5$ depends only on $g_0$.

On the other hand, since $n + \Delta_{g_0} u > 0$ on $M$, applying the Green function of $g_0$, we can get
\begin{equation}
\label{c0-estimate-7}
\sup_M u \,\le\, \frac{1}{V} \int_M u \omega_{g_0}^n + C_6,
\end{equation}
where $V=\int_M \omega_{g_0}^n$ and $C_6$ depends only on $g_0$.

By our assumption on $u$, we get from the above
$$\int_M u \omega_{g_0}^n \,\le\, V\, (1 + C_6).$$
Combining this with \eqref{c0-estimate-5} and \eqref{c0-estimate-6}, we obtain an a prior estimate on
$||u||_{C^0}$ and the proposition is proved in the case that $\sup_M u =-1$.

In general, if $u$ is a solution of \eqref{continuity}, then $\bar u := u - \sup_M u -1$ is also a solution.
Applying the above discussion, we have
$$\sup_M u - \inf_M u \le C_7,$$
where $C_7$ is a uniform constant. Therefore, if $u$ satisfies $\int_M u\, \omega_{g_0}^n =0$, then by \eqref{c0-estimate-7},
$$\sup_M u \le C_6.$$
Hence, we have a uniform estimate on $||u||_{C^0}$ as required by the proposition.
\end{proof}

\subsection{The higher order estimates}
In order to establish the existence of V-solitons, we need higher order estimates for solutions of \eqref{continuity}.
Based on the known theory on the $C^{2,\alpha}$-estimate for complex Monge-Ampere equations, we only need an a prior $C^2$-estimate.

The following is trivial.
\begin{lemma}\label{uniformhessian} Let $u$ be a solution of \eqref{continuity}, then
$||\p\bar \p u||\,\leq \,\max \{ n+\Delta_{g_0}  u, n\}$
\end{lemma}

Therefore, in order to derive an a priori $C^2$-estimate, we only need to have a
$C^0$-estimate for $\Delta _{g_0} u $.
This is similar to the second-order estimate in the proof for the Calabi-Yau theorem.
However, because of the extra term involving
$|V|_u^2$, the proof in our case is much more tricky and lengthy. This will be in our forthcoming paper.

\subsection{Uniqueness of scalar V-soliton equation}
Using the Maximum principle, one can easily show
the following:

\begin{theorem}
\label{uniqueness}
Let $(M, g_0)$ be a compact K\"ahler manifold with boundary $\partial M$ and a $S^1$-symmetry induced by
a Hamiltonian field $V$. Then there is at most one solution of \eqref{twistedmongeampere} with given boundary value,
namely, if $u_1$ and $u_2$ are $S^1$-invariant solutions of \eqref{twistedmongeampere} with $u_1=u_2$ along $\partial M$, then
$u_1\equiv u_2$ on $M$.
\end{theorem}

\section{Further directions}
In this section, we discuss possible applications of our new correspondence and some further research problems.

\subsection{Boundary value problem for V-soliton metrics}
First we certainly concern the existence problem of V-soliton metrics. This is amount to solving the scalar V-soliton equation
\eqref{twistedmongeampere}. We expect: {\it Given a complete K\"ahler manifold $(M,g_0)$ with
boundary $\partial M$ and finite geometry at $\infty$. Suppose that it admits a $S^1$-symmetry generated by a Hamiltonian field $V$,
then for any reasonably ''nice'' boundary value $\varphi$ along $\partial M$, there is a unique solution $u$ of \eqref{twistedmongeampere}
on $M$ such that $u|_{\partial M} = \varphi$}.

In \cite{gatian}, we will provide a solution to this existence problem in the case that $M$ is compact or an ALE space.
The solution we obtain will be in $C^{1,1}$ in general, but it should be smooth outside the zero set of $V$ as an application of the
known regularity theory for complex Monge-Ampere equations. It will be a more challenging problem to study the regularity of such a solution
near the zero set of $V$.

\subsection{Finite-time singularities of the K\"ahler-Ricci flow}\label{subsection:finitetime}
Our new correspondence may be applied to studying singularity formation of the K\"ahler-Ricci flow: Let $(M,g)$ be a K\"ahler
manifold with a $S^1$-symmetry generated by a Hamiltonian field $V$. We further assume that $g$ is a V-soliton metric (i.e., it satisfies eq. \eqref{eq:twist-einstein2}).
Let $\mu: M\mapsto \RR$ be the associated moment map, i.e., the Hamiltonian function of $V$. Put $Cr(\mu)$ to be the set of critical
values of $\mu$. Then $\RR\backslash Cr(\mu)$ is a disjoint union of consecutive open intervals $I_a$ ($a\in \ZZ$). For each interval
$I_a$, symplectic quotients $X_\tau$ for $\tau \in I_a$ are the same complex manifold, but $X_\tau$ changes when $\tau$ crosses critical values in $Cr(\mu)$. Usually, $X_\tau$ and $X_{\tau'}$ are related to each other by so called flips when $\tau$ and $\tau'$ are in two different, but
consecutive, intervals. By studying how $g$ descends to $X_\tau$ and $X_{\tau'}$, we can analyze how the K\"ahler-Ricci flow transforms under
flips. Let us illustrate this by means of an example.

Let $\CC^*$ act on $M:=\CC^{l+m}$ by:
$$ t(z_1,\ldots,z_{l+m})=(t^{a_1}\cdot z_1,\ldots,t^{a_l} \cdot
z_l,t^{-a_{l+1}}\cdot z_{l+1},\ldots, t^{-a_{l+m}}\cdot z_{l+m}),$$
\noindent
where $a_1,\ldots,a_l,a_{l+1},\ldots,a_{l+m}>0$ are positive integers.
This action is Hamiltonian with respect to the standard K\"ahler structure on $\CC^{l+m}$ with the Hamiltonian $\mu : \CC ^{l+m} \to \RR$:
$$\mu (z_1,\cdots, z_{l+m})= \sum _{i=1}^{l}a_i |z_i|^2 - \sum _{i=l+1}^{l+m}a_i |z_i|^2$$
One can easily see that $\tau=0$ is the only critical value. Therefore,
the symplectic quotients $X_{\tau}:=\mu ^{-1}(\tau)/S^1$ are all isomorphic to a fixed variety $X^-$ for $\tau <0$ and to a variety $X^+$ for $\tau >0$.
Furthermore, the natural bi-rational map $\phi : X^- \dashrightarrow X^+$
is a flip for $l, m\geq 2$ replacing via surgery a neighborhood of ${\PP}^{l-1}\subset X^- $
with a neighborhood of ${\PP}^{m-1}\subset X^+ $.
For $l=1, m\geq 2$ and $a_1=\cdots=a_{l+m}=1$, $\phi$ is a blow-down, and for
$l\geq 2, m=1$ and $a_1=\cdots=a_{l+m}=1$ it is a blow-up. For $l\geq 2$ and $a_1=\cdots=a_{l+m}=1$ it is a flip or flop (e.g., $l=m=2$ gives rise to a flop).
Another important special case is when $l=m=2$, $a_1=2$ and $a_2=\cdots =a_4=1$: This is the first non-trivial flip in the Francia series.

Let us consider the simplest case in this context: $l=1$ and  $a_1=\cdots=a_{l+m}=1$. For $\tau<0$, $X_\tau$ is the $S^1$-quotient of
$$\{(z_1,z_2,\cdots,z_{m+1})\,|\, |z_1|^2 + |\tau| = \sum _{i=2}^{m+1} |z_i|^2 \}$$
which is the blow-up of $\CC ^m$ at $(0,0)$. On the other hand, one can see easily that $X_\tau = \CC^m$
for $\tau >0$. Let us find a special $V$-soliton metric $g$ on $\CC^{n+1}$ of the form
$$\omega_g = \frac{\sqrt{-1}}{2} \p\bar \p u,~~~~u=  (|z_2|^2 +\cdots + |z_{m+1}|^2)\,h(|z_1|^2).$$
The holomorphic field whose imaginary part equals $V$ is given by
$$W\,=\, - z_1\frac{\p}{\p z_1} + z_2 \frac{\p}{\p z_2} + \cdots + z_{m+1} \frac{\p}{\p z_{m+1}}.$$
Then the scalar $V$-soliton equation is equivalent to the following:
$$\det (u_{i\bar j}) \,=\, u_{1\bar 1} |z_1|^2  - \sum_{j=2}^{m+1} (u_{1\bar j} z_1 \bar z_j + u_{j \bar 1} z_j\bar z_1)
+ \sum_{i,j\ge 2} u_{i\bar j}z_i\bar z_j,$$
where $u_{i\bar j}= \frac{\p^2 u}{\p z_i\p \bar z_j}$.

If we set $r:=|z_1|^2$ and $\rho=|z_2|^2 +\cdots + |z_{m+1}|^2$,
one can verify directly that $u= \rho h(r)$ satisfies the V-soliton equation
if $h$ satisfies the following ODE:
$$\left (  r h h''  - r (h')^2  + h h' \right ) h^{m-1} = r^2  h'' - r h' + h,$$
where $r\in [0,\infty)$ and $h$ is a function of $r$. Given a solution $h$ of this equation, we obtain a $V$-soliton metric $g$
on $\CC^{m+1}$ which descends to a solution of the K\"ahler-Ricci flow for $\tau < 0$ and converges to
a smooth K\"ahler metric $g_0$  on $\CC^m\backslash \{0\}$. One can show easily that $\CC^m$ is the metric completion of
$\CC^n\backslash \{0\}$ by $g_0$. This can be used to verify the first conjecture on finite-time singularity for the K\"ahler-Ricci flow
in \cite{ti4} and \cite{songtian} in the special case of a blow-up of a smooth manifold. 
In fact, in order to see how the K\"ahler-Ricci flow behaves under
blowing-down of a $\CC P^{m-1}$, it suffices to find a solution of the above ODE near $r=0$. By using the power series method, one can find a
local solution $h$ of the above ODE starting with
$$h (r) = 1 + r - \frac{m-1}{4} r^2 + O(r^3).$$
One can use this explicit solution to see that $\CC^m$ is the metric completion of
$\CC^m\backslash \{0\}$ by $g_0$ as claimed above. Similarly, one can find special solutions for the V-soliton equation
in the general case that both $l, m > 1$.

In fact, this example in the case of $m=2$ provides the basic picture of finite-time singularity
for the K\"ahler-Ricci flow on complex surfaces. Let us elaborate more on this: Let $X$ be a complex surface
and $g(0)$ be a K\"ahler metric. Then there is a unique solution $g(t)$ of the K\"ahler-Ricci flow on $[0, T)$,
where $T$ is either $\infty$ or the first time when $[\omega_{g(0)} ] - t c_1(X)$ fails to be positive. If $T< \infty$
and $([\omega_{g(0)}] - T c_1(X))^2 > 0$, then as $t\to T$, $g(t)$ converges to $g(T)$ outside finitely many disjoint
rational curves $C_1,\cdots, C_k$ of self-intersection number $-1$. For simplicity, assume that $k=1$.
By blowing down $C_1$, we get a new complex surface
$\bar X$ with $p$ corresponding to the blow-down $C_1$. We can also extend $g(T)$ to be a solution $g(t)$ of
the K\"ahler-Ricci flow on $\bar X$ for $t\in [T, T+\epsilon]$ for some $\epsilon >0$.
Let $U$ be a small neighborhood of $p$ and $\tilde U$ be the blow-up of $U$ at $p$. Then
$\tilde U$ to be a neighborhood of $C_1$ in $X$, moreover, we can identify
$\tilde U\times (T-\epsilon, T)\cup U\times [T, T+\epsilon)$ as a quotient of a neighborhood $W\subset \CC^3$ of $0$ by the $S^1$-action.
The solution $g(t)$ lifts to a V-soliton metric $\bar g$ on $W\backslash \{0\}$. By solving the boundary value problem
for the scalar V-soliton equation on $W$, we should be able to extend $\bar g$ on $W$. Then by studying how $\bar g$ descends to $U$,
we may prove that $\bar X$ is the metric completion of $\bar X\backslash \{p\}$ by $g(T)$. This verifies the first conjecture
on finite-time singularity for the K\"ahler-Ricci flow
in \cite{ti4} and \cite{songtian} for complex surfaces. 
Of course, the above discussion just provides a plausible approach. Details remain to be checked.

We believe that this actually provides an effective approach to studying finite-time singularity of the K\"ahler-Ricci flow in all dimensions,
at least for all those flips which can be achieved through variations of symplectic quotients. Indeed, many flips can be achieved in this way.
This allows us to carry out a geometric Minimal Model Program using the Ricci flow with "surgeries". The first step in this program is to understand
solutions for the $V$-soliton equation on a manifold with boundary. This will be the subject of \cite{gatian}.

\subsection{K\"ahler-Ricci flow on Fano manifolds}
Another possible application of the $V$-soliton equation is to study the K\"ahler-Ricci flow on Fano manifolds. Let $X$ be a Fano manifold
and $g_0$ be a K\"ahler metric with its K\"ahler class equal to $c_1(M)$. It is known that the normalized K\"ahler-Ricci flow
$$\frac{\p g}{\p t} \,=\, - {\rm Ric}(g) + g,~~~g(0)=g_0$$
has a global solution $g(t)$ for all $t>0$. A long-standing problem is on the convergence of $g(t)$ as $t$ goes to $\infty$. The folklore conjecture
is that $g(t)$ converges to a K\"ahler-Ricci soliton (possibly with mild singularity along a subvariety of complex codimension at least $2$).
Our new correspondence may provide a method of proving this conjecture. By Theorem \ref{theorem:converse}, there is a K\"ahler metric
$\bar g(\cdot, z)$ on $M=X\times \{z\in \CC\,|\, |z|\ge 1\}$ satisfying:

\noindent
(1) $\tau = e^t - 1$;

\noindent
(2) $\bar g$ is invariant under the standard $S^1$-action of $\CC$ by rotations;

\noindent
(3) $g(t)$ is the symplectic quotient of $\bar g$ on $X$.

\noindent
(4) $\bar g$ satisfies the $V$-soliton equation:
$${\rm Ric }(\bar g) +\frac{\sqrt{-1}}{2} \p\bar \p (- \log w + f) = \omega_{\bar g}, $$
where $w$ is the inverse of the squared norm of $V$ ($w=(|V|_g^2)^{-1}$) given on $X\times [0,\infty)$ by
$$w=4^{-1}\,e^{- 2 t} \left (R(g(t)) - n + \frac{\p  f(e^t-1)}{\p t} \right )\, >\, 0.$$
This is the same as
$$4^{-1}\,(1+\tau)^{-2} \,\left (R(g(t)) - n + (1+\tau) \frac{\p  f(\tau)}{\p \tau} \right )\, >\, 0.$$
That it is possible to find such an $f$ is insured by a result of Perelman's  to the effect that $R(g(t))$ is bounded (cf. \cite{sesumtian}).
Such an $f$ is not unique, so we may choose
one that is more convenient to us. For instance, if
$c$ is the lower bound of $R(g(t))$, we choose $f=(n-c) \log (1+\tau) + 2 (1+\tau)^2$. Then
$$w \,=\, 4^{-1} (1+\tau)^{-2} \,\left (R(g(t)) - n + \frac{\p  f(e^t-1)}{\p t}\right )\, \sim\,  1.$$
It follows that at $\infty$ of $\CC$, in polar coordinates $z=(\tau,\varphi)$, we have
$$g\,\sim\,\lim_{t\to \infty} g(t) + d\tau^2 +  d\varphi^2.$$
If $g(t)$ converges to a K\"ahler-Einstein metric $g_{KE}$ as $t$ tends to $\infty$, then
$\bar g$ can be extended across $X\times \{\infty\}$ by adding $g_{KE}$. Or equivalently, given any
sequence $\{t_i\}$ with $\lim t_i = \infty$, then $(X\times \CC, \bar g (t+ t_i))$ converges to
the product of $g_{KE}$ with the cylinder metric $d\tau^2 + d\varphi^2$ as $t_i$ tends to $\infty$.

In general, it is plausible that the above chosen $\bar g$ can be extended across $\infty$ of $\CC$ modulo a family of diffeomorphisms of $X$
or equivalently, $(X\times \CC, \bar g (t+ t_i))$ converges to
the product of a limiting metric with the cylinder metric $d\tau^2 + d\varphi^2$.
Therefore, the above folklore conjecture is closely related to how $(X\times \CC, \bar g)$ behaves at the $\infty$ of $\CC$ and whether or
not it can be compactified. We conjecture that {\it $(X\times \CC, \bar g (t+ t_i))$ converges to
the product of a K\"ahler-Ricci soliton with $d\tau^2 + d\varphi^2$ modulo diffeomorphisms.}
Based on this idea and assuming the analyticity, Arezzo and La Nave (cf. \cite{gabcla}) studied the case that the central fiber of a (non-trivial) special degeneration $\cX \to \Delta$ admits a K\"ahler-Einstein metric. In a forthcoming paper (cf. \cite{aregatian}), we will discuss this in more details.

\subsection{$V$-solitons and geodesics in the space of K\"ahler metrics}

On a K\"ahler manifold $X$, each $(1,1)$-form cohomologous to $\omega$ takes the form $\omega + \sqrt{-1}
\partial\bar\partial f$ for some $f\in C^{\infty}(X)$. Therefore, the space of
all K\"ahler metrics in the class $[\omega]$ can be identified with
\begin{equation}\label{twistedmongeampere}
\mathcal{H}=\{\phi\in
C^\infty(X)~|~\omega +\frac{\sqrt{-1}}{2} \p \bar \p \phi >0\}/\thicksim.
\end{equation}
where $\phi \thicksim \phi'$ if and only if they are different by addition of a constant.

Given $\phi\in \mathcal{H}$,  the formal
tangent space
$$T_{\phi}\mathcal{H}\,=\, \{\psi \in C^\infty(X)_0~|~\int_X\psi\omega_\phi^n=0\},$$
where $\omega _{\phi}:=\omega+\frac{\sqrt{-1}}{2} \partial\bar\partial\phi.$
There is an natural metric, introduced by T. Mabuchi, on $\mathcal{H}$ as follows: Let $\psi_1,\psi_2\in T_\phi\mathcal{H}$, define
\begin{equation}
<\psi_1,\psi_2>_{\phi}=\frac{1}{n!}\,\int_X \psi_1\psi_2\,\omega_\phi^n.
\end{equation}
Given a smooth curve $\phi(t):[a,b]\mapsto \mathcal{H}$, set $\phi(x,t):= \phi(t)(x)$. This can be considered as
a function on $X\times S^1\times [a,b]$ which is $S^1$-invariant.
Then the geodesic equation for the above $L^2$ metric is equivalent to the following Homogeneous
complex Monge-Ampere (in short HCMA) equation on $M= X\times  S^1\times[a,b]$ (cf. \cite{HCMAsemmes}, \cite{donaldsonholdisk} and \cite{ChenTian}):

\begin{equation}\label{hcma}
(\omega +\frac{\sqrt{-1}}{2}  \p_M \bar \p _M \phi )^{n+1}=0.
\end{equation}
On the other hand, consider the K\"ahler-Ricci flow on the K\"ahler manifold $X\times [0,\infty)$:
$$\frac{\p g}{\p t} \,=\, - {\rm Ric}(g) +  g,~~~g(0)=g_0.$$
By Theorem \ref{theorem:converse}, there is a lifting metric $\bar g(\cdot, \tau)$ on $M=X\times S^1\times [0,\infty)$ satisfying the following

\noindent
(1) $\tau = e^{ t} - 1$;

\noindent
(2) $\bar g$ is invariant under the obvious $S^1$-action;

\noindent
(3) $g(t)$ is the symplectic quotient of $\bar g$ on $X$.

\noindent
(4) For some $f$, $\bar g$ satisfies the $V$-soliton equation:
$${\rm Ric }(\bar g) +\frac{\sqrt{-1}}{2} \p\bar \p ( \log |V|_g^2 + f) = \omega_{\bar g}, $$
where $V$ is the vector d=field generating the $S^1$-action.
As we have shown in the above, the $V$-soliton equation can be reduced to
the scalar $V$-soliton equation:
\begin{equation}\label{twistedmongeampere}
(\omega_{g_0} +\frac{\sqrt{-1}}{2} \p _M \bar \p _M \phi )^n = |V| _{g}^2
e^{F+f-\phi}\,\omega_{g_0}^n,
\end{equation}
where $F$ is given by
$${\rm Ric}(g_0) -\lambda \omega_{g_0} = \frac{\sqrt{-1}}{2}\p _M  \bar \p _M F.$$
Since $V$ tends to $0$ as $\tau$ goes to $\infty$, we may expect that the solutions to \eqref{eq:v-twisted-MA} are asymptotic to the solutions of  equation \eqref{hcma}, more precisely, we expect that the solution $g(t)$ of the Ricci flow is asymptotic to a geodesic ray. This can be a future research topic.

\end{document}